\documentclass{amsart}
\advance\textwidth 1cm
\advance\oddsidemargin -0.5cm
\advance\textheight 0.5cm
\advance\topmargin -0.3cm
\title{Solovay's Relative Consistency Proof for FIM and BI}
\author{Joan Rand Moschovakis  \\ Occidental College (Emerita) \\ \email{joan.rand@gmail.com}}

\usepackage{amsmath}

\begin{document}
\maketitle

Robert M. Solovay is a classical mathematician and logician whose many contributions to set theory, proof theory, cryptography and computer science are well known.  Although he ``thinks classically,'' occasionally a question involving intuitionistic logic attracts his attention.  This historical note documents how he answered one such question and contributed a rather surprising result to the metamathematics of intuitionistic analysis.

In 2002 Solovay wondered about the consistency strength of Kleene's system {\bf I} of intuitionistic analysis (cf. \cite{klve1965}) (which he called ``{\bf FIM}''), in relation to classical systems such as {\bf ATR} and {\bf $\Pi^1_1$-CA} studied by Simpson in \cite{si2009}.  Kleene had proved in \cite{kl1969} that {\bf I} is consistent relative to its neutral, classically correct subsystem {\bf B}, which Solovay called ``{\bf BSK}''.  In this exposition Solovay's acronyms for {\bf I} and {\bf B} will replace Kleene's from now on.

Solovay conjectured that {\bf FIM} has the same consistency strength as a classical theory {\bf BI} (``Bar-Induction'') he described as follows in an email message to the author on 30 July 2002:\footnote{The quotations are verbatim, except for correcting trivial typographical errors (e.g. ``consistenciency''), rendering sub- and superscripts in latex, and then replacing e.g. ``Pi$^1_1$'' by ``$\Pi^1_1$.''}

\begin{quotation}

The logic is classical.  There are two sorts of variables: lower case letters stand for the number variables;  upper case letters for the set variables.

There is a binary predicate = in two flavors: for equality of numbers or equality of sets; there is a binary epsilon relation; there are the usual function symbols from Peano: 0, S, +, $*$.

There are the usual Peano axioms.  Induction is in the strong form. [Arbitrary formulas of the language are allowed.]  One has extensionality for sets.

One has arithmetic comprehension.  The set of n such that Phi(n) exists [for each Phi without bound set variables.]

The key axiom asserts that if R is a binary relation R which is a linear ordering and has the property that for every non-empty subset of its field there is an R-least element, then one has *full* R-induction.  [*full* means that for any formula Phi expressible in the language, the corresponding ``least element with respect to R'' principle holds.]

That completes my description of Bar-Induction.

These facts are known about the consistency strength of Bar-Induction:

1) It proves the consistency of ATR$_0$.  [Roughly ATR$_0$ says that if R is a well-ordering of omega, and X is a set, then one can define the hyperarithmetical hierarchy relativized to X with indices chosen from R.]

2)  $\Pi_1^1$-Comprehension proves that Bar-Induction has a beta model [one in which all the things that the model believes are well-orderings are indeed well-orderings].

It looks quite routine to me to carry out the usual realizability of FIM within BI [henceforth, my abbreviation for Bar-Induction].  This will give a relative consistency proof for FIM relative to BI.

To go the converse direction the obvious try is the ``negative interpretation'' of Godel.  But I seem to get something like BI is consistent if BSK [ = FIM - Continuity] + ``Markov's principle'' is consistent.

So one needs BSK + Markov's principle is consistent relative to BSK. \ldots

\end{quotation}

For comparison with {\bf BI} a brief description of {\bf FIM} is needed here.  The logic is intuitionistic.  There are two sorts of variables: lower case Latin letters $\mathrm{a, b, c, \ldots, x, y, z, a_1,\ldots}$ stand for the number variables; lower case Greek letters $\mathrm{\alpha, \beta, \ldots}$ stand for the (one-place number-theoretic) function variables (the ``choice sequence'' variables, in intuitionistic terminology).

There are the usual function symbols 0, S, +, $\cdot$, plus a finite number of symbols for additional primitive recursive functions and functionals sufficient to formalize the theory of recursive partial functionals (as Kleene did in \cite{kl1969}).  There is a binary predicate = for equality of numbers; if $\mathrm{s,t}$ are terms (of type 0) then $\mathrm{s = t}$ is a prime formula. Parentheses indicate application, e.g. $\mathrm{\alpha(x)}$ is a term.  Equality of functions is defined extensionally by $\mathrm{(\alpha = \beta) \equiv \forall x (\alpha(x) = \beta(x))}$.  Church's $\lambda$, which can be eliminated, allows constructing functors (terms of type 1) from terms, e.g. $\mathrm{\lambda x. 0}$ represents the identically zero function, and there is a $\lambda$-reduction schema.

There are the usual Peano axioms and the defining axioms for the additional function(al) symbols.  Induction is in the strong form (arbitrary formulas of the language are allowed). There is an open equality axiom $\mathrm{x = y \rightarrow \alpha(x) = \alpha(y)}$, and equality axioms for the function(al) constants are provable in the subsystem {\bf IA$_1$} (``two-sorted intuitionistic arithmetic'') of {\bf FIM} outlined so far.

Now Kleene did not restrict countable choice to arithmetical formulas. Full countable choice for numbers AC$_{00}$:
\[\mathrm{\forall x \exists y A(x,y) \rightarrow \exists \beta \forall x A(x,\beta(x))}\]
is a theorem ($^*$2.2 in \cite{klve1965}) of {\bf BSK} = {\bf IA$_1$} + AC$_{01}$ + BI!, where AC$_{01}$ ($^\mathrm{x}$2.1 in \cite{klve1965}) is the axiom schema of countable choice for functions:
\[\mathrm{\forall x \exists \alpha A(x,\alpha) \rightarrow \exists \beta \forall x A(x,\lambda y. \beta(2^x \cdot 3^y))}\]
and BI! is an axiom schema ($^\mathrm{x}$26.3c in \cite{klve1965}) of bar induction, to be described later.

In \cite{kl1969} Kleene formalized the theory of recursive partial functionals in {\bf IA$_1$} plus the countable comprehension schema AC$_{00}$!:
\[\mathrm{\forall x \exists ! y A(x,y) \rightarrow \exists \beta \forall x A(x,\beta(x))},\]
where in general $\mathrm{\exists ! y B(y) \equiv \exists y B(y) ~\&~ \forall x \forall y (B(x) ~\&~ B(y) \rightarrow x = y)}$.
Anne Troelstra observed (p. 73 of \cite{tr1973}) that quantifier-free countable choice qf-AC$_{00}$ (AC$_{00}$ for formulas $\mathrm{A(x,y)}$ containing only bounded numerical quantifiers, with parameters allowed) suffices for this purpose.  Like Troelstra's system {\bf EL}, the subsystem {\bf IRA} = {\bf IA$_1$} + qf-AC$_{00}$ of {\bf BSK} precisely expresses intuitionistic recursive analysis (cf. Vafeiadou \cite{gvfms}).

Kleene's formalization of function-realizability culminated in Theorem 50 of \cite{kl1969}, which had the simple consistency of {\bf FIM} relative to {\bf IA$_1$} + AC$_{00}$! + BI! as a corollary. Troelstra's remark on p. 208 of \cite{tr1973} suggested that AC$_{00}$! might be weakened to qf-AC$_{00}$ here too.  Since arithmetic comprehension is expressible in the language of {\bf FIM} by restricting AC$_{00}$! to arithmetic formulas $\mathrm{A(x,y)}$, and quantifier-free formulas are arithmetic, Solovay's claim that {\bf FIM} is consistent relative to {\bf BI} was not very surprising.

It did not seem at all obvious, however, that arithmetic comprehension could be negatively interpreted in {\bf BSK} + MP.  Working out the details required choosing correctly among Kleene's four versions (26.3a-d in \cite{klve1965}, described following the quotation) of the bar-induction axiom, as Solovay observed in an email message of 1 August 2002:

\begin{quotation}

BI is the system referred to by Simpson as $\Pi^1_\infty$-TI$_0$.

Now the various forms of 26.3 are not equivalent in BI.  The reason is that we don't have full comprehension or x2.2, so it's hard to get a real as in 26.3b from a decidable predicate.  So we have 26.3b in BI but *not* 26.3a.

Put it another way:

In BI a linear-ordering coded by a real that happens to be a well-ordering one can do induction along.  But if the linear ordering is just given by a formula of our language, BI says nothing.  [Of course, in the presence of full comprehension the distinction between ``coded by a real'' and ``given by a formula'' evaporates.]

\end{quotation}

Letting $\mathrm{w}$ vary over sequence numbers (primitive recursive codes for finite sequences with a primitive recursive length function $\mathrm{lh(w)}$, where 1 codes the empty sequence, $\mathrm{2^{n+1}}$ codes the sequence $(n)$,  and $\mathrm{*}$ denotes concatenation), Kleene's $^\mathrm{x}$26.3a can be compactly stated as follows:
\[\mathrm{\forall w (R(w) \vee \neg R(w)) ~\&~ \forall \alpha \exists x R(\overline{\alpha}(x)) ~\&~ \forall w (R(w) \rightarrow A(w)) ~\&~ \forall w (\forall n A(w * 2^{n+1}) \rightarrow A(w)) \rightarrow A(1)},\]
while $^\mathrm{x}$26.3b (henceforth also referred to as ``BI$_1$'') can be abbreviated by
\[\mathrm{\forall \alpha \exists x \rho(\overline{\alpha}(x)) = 0 ~\&~ \forall w (\rho(w) = 0 \rightarrow A(w)) ~\&~ \forall w (\forall n A(w * 2^{n+1}) \rightarrow A(w)) \rightarrow A(1)},\]
where $\mathrm{\rho}$ represents an arbitrary element of Baire space (``a real'').\footnote{Kleene considered two other forms of bar induction, $^\mathrm{x}$26.3c or BI! (like $^\mathrm{x}$26.3a but replacing the hypotheses $\mathrm{\forall w (R(w) \vee \neg R(w)) ~\&~ \forall \alpha \exists x R(\overline{\alpha}(x))}$ by $\mathrm{\forall \alpha \exists ! x R(\overline{\alpha}(x))}$).  BI$_1$ is weaker than the other forms over {\bf IRA}, but all four are equivalent over {\bf IA$_1$} + AC$_{00}$!.}

The version of Markov's Principle referred to here as ``MP'' is the strong analytical form
\[\mathrm{\forall \alpha [\neg \forall \alpha \neg (\alpha(x) = 0) \rightarrow \exists x \alpha(x) = 0]}.\]
Kleene's formalization of function-realizability in \cite{kl1969} established that {\bf FIM} + MP is simply consistent relative to a subsystem of {\bf BSK} + MP with countable choice weakened to countable comprehension.  Solovay observed that arithmetical comprehension is enough.

In an email message of 20 August 2002 he outlined a proof in primitive recursive arithmetic {\bf PRA} that {\bf BSK} + MP is consistent relative to {\bf BI}.

\begin{quotation}

I've thought through about 95$\%$ of the details on the following:

PRA proves [Con(BI) implies Con(BSK + MP)].

The proof basically is just carrying out the function realizability for BSK inside of BI.  One actually proves the following in PRA.

There is a primitive recursive function which takes as input a proof of a sentence, A, in BSK + MP and outputs the following:

(a) a Godel number e of a partial recursive function [mapping numbers to numbers];

(b) a proof in BI that the function, call it $\phi_n$, is total.

(c) a proof in BI that $\phi_n$ realizes A.

Of course, much less would suffice for the relative consistency proof.  We could get by with the output

(d) A proof in BI that there is some function, $\alpha$ [not necessarily recursive] that realizes A.

I remark that with realizability as presented in Kleene one can not sharpen (a) and (b) to

(e) a ``primitive recursive Godel number'' of some primitive recursive function from numbers to numbers.  [The reason is that if we could get (e) we could find a primitive recursive algorithm for every provably recursive function of the theory BI, and Ackerman's fcn. is such a provably recursive function.]

\end{quotation}

To complete his theorem Solovay needed to prove {\bf BI} consistent relative to {\bf BSK} + MP.   While BI$_1$ (Kleene's $^\mathrm{x}$26.3b) is obviously negatively interpretable in  {\bf BSK} + MP, and most other axioms of {\bf BI} are negatively interpretable in {\bf BSK}, arithmetical comprehension presented a problem. After one false start, on 31 August 2002 Solovay sent the author the following message in which the role of the {\em intended good classical models} is explicitly acknowledged:

\begin{quotation}

I think I've got it right.  Recall that the plan is to define a classical theory BI- such that

(1) BI- has the same theorems as BI and this is finitistically provable.

(2) BI- is designed to be easy to negatively interpret in BSK + MP.

I have a new candidate for BI- for which I have checked (1).  I'm optimistic re (2) but I have not yet thought through the details.

I will consider a variant of BI where the type 1 variables range over functions from omega to omega.

Without describing BI with utter precision, it includes:

(1) The usual axioms and operations of Peano Arithmetic;

(2) Full induction for any formula expressible in the language [parameters allowed].

(3) Axioms that insure that the type 1 variables are closed under Turing jump and ``recursive in'' and contain all recursive functions.

(4) The key axiom of Bar-Induction:  If R is a linear ordering on omega, and there is no descending chain through R given by one of the type 1 functions, then one has induction over R for arbitrary formulas of the language [parameters allowed].

BI- is obtained from BI by replacing (3) and (4) by suitable variants:

(3$'$) We require that the type 1 functions contain all primitive recursive functions and that if $\alpha$ and $\beta$ are type 1 functions and that $\gamma$ is primitive recursive in $\alpha$ and $\beta$ then $\gamma$ is a type 1 function.

[Of course, I'm being sloppy here and implicitly describing axioms by describing what the intended good models of the theory are.]

(4$'$) Axiom x26.3b of Kleene's FIM. [Caution: In the current classical context, it makes quite a difference which version of 26.3 one takes.]

\ldots

As I said, I'm claiming that BI- and BI have the same theorems and that this can be proved in PRA.

\end{quotation}

Now {\bf BI-} looks like a classical version of a subsystem of {\bf IRA} + BI$_1$.  In response to a request for a proof in {\bf BSK} + MP of the negative interpretation of arithmetic comprehension, Solovay responded on 2 September 2002 with a beautiful argument for the result referred to as ``Solovay's Lemma'' in \cite{jrm2003}.\footnote{The author has since used variations of his argument for several purposes in the context of intuitionistic analysis, but Solovay has not published it himself. This is a LaTeX transcription of his original proof.}

\begin{quotation}

I haven't tried for a direct proof.  But perhaps the place to start is analysing the proof of arithmetic comprehension in BI-.

Here is a sketch of my argument.

Let $\mathrm{\alpha: \omega \rightarrow \omega}$.  We aim to prove the existence of a $\mathrm{\beta}$ with the following properties:

1) $\mathrm{\beta(2n) = \alpha(n)}$;

2) $\mathrm{\beta(2n+1) > 0}$ iff $\mathrm{\exists y T^\alpha(n, n, y)}$.

3) If $\mathrm{\beta(2n+1) > 0}$ then it equals $\mathrm{y + 1}$ where $\mathrm{y}$ is least such that $\mathrm{T^\alpha(n,n,y)}$.

We first define the following $\mathrm{\rho}$ which will be uniformly primitive recursive in $\mathrm{\alpha}$:

1) If $\mathrm{s}$ is not a sequence number then $\mathrm{\rho(s) = 1}$.

2) Now let $\mathrm{s}$ be a sequence number.  If for some $\mathrm{j < lh(s)}$, we have $\mathrm{j = 2k}$ and $\mathrm{(s)_j \neq \alpha(k)}$, then $\mathrm{\rho(s) = 0}$;

OR 3) if for some $\mathrm{j < lh(s)}$ we have $\mathrm{j = 2k+1}$ and $\mathrm{(s)_j = 0}$ and $\mathrm{\exists y \leq lh(s) T^\alpha(k, k, y)}$ then $\mathrm{\rho(s) = 0}$:

OR 4) if for some $\mathrm{j < lh(s)}$, we have $\mathrm{j = 2k+1, (s)_j = m + 1}$ and $\mathrm{m}$ is not the least $\mathrm{y}$ such that $\mathrm{T^\alpha(k,k,y)}$ then $\mathrm{\rho(s) = 0}$.

OTHERWISE $\mathrm{\rho(s) = 1}$.

Now I describe the predicate $\mathrm{A(x)}$.  [For use in the Brouwer principle axiom.]

not $\mathrm{A(s)}$ iff

1) $\mathrm{s}$ is a sequence number;

2) let $\mathrm{j = 2k < lh(s)}$.  Then $\mathrm{(s)_j = \alpha(k)}$.

3) let $\mathrm{j = 2k+1, \; j < lh(s)}$.  Then $\mathrm{(s)_j > 0}$ iff $\mathrm{\exists y T^\alpha(k,k,y)}$.  If so letting $\mathrm{y_k}$ be the least such $\mathrm{y}$ we have $\mathrm{(s)_j = y_k+1}$.

From the fact that not $\mathrm{A(1)}$ we conclude by bar induction that $\mathrm{\exists \gamma\forall n \rho(\overline{\gamma}(n)) > 0}$. [Recall that we are reasoning in the ``classical'' system BI-.] But then it is easy to see that this $\mathrm{\gamma}$ is our desired $\mathrm{\beta}$.

\end{quotation}

This proof, which guarantees that the range of the type 1 variables in an omega-model of {\bf BI-} is closed under the Turing jump, uses only primitive recursive comprehension and classical bar induction. Arithmetic comprehension follows easily by formula induction.

To complete his proof that Con({\bf BSK} + MP) implies Con({\bf BI}) Solovay confirmed in another email message on 2 September 2002 that {\bf BI-} can be negatively interpreted in {\bf BSK} + MP:

\begin{quotation}

The theory BI- doesn't explicitly have arithmetic comprehension among the axioms.  Instead it's a theorem.  But the set of things whose negative interpretation is a theorem of BSK + MP is closed under classical deducibility.  So it is enough to check the axioms of BI- have negative interpretations that are theorems of BSK + MP.

The only problematical axioms are the instances of the version of Brouwer's principle that are axioms of BI-.  Such an axiom has four clauses:

Hypotheses:

1) Every $\mathrm{\alpha}$ hits a bar.

2) $\mathrm{A}$ holds at bars.

3) $\mathrm{A}$ holds at $\mathrm{s}$ if it holds at all one step extensions of $\mathrm{s}$.

Conclusion:

4) $\mathrm{A}$ holds at 1.

Now the negative interpretations of clauses 2) through 4) transform into things of the same shape. [The $\mathrm{A}$ gets replaced by its negative interpretation.]

But 1) gives trouble.  The inner existential number quantifier gets negatively replaced.  But MP allows us to restore this existential number quantifier.  This is crucial since the version of Brouwer available in BSK has the existential number quantifier and not its negative transform.

So roughly the negative transform puts a crimp in clause 1) and MP allows us to remove it.

\end{quotation}

The result is an elegant, elementary proof of the relative consistency of formal systems of classical and intuitionistic analysis with bar induction (and Markov's Principle):

\vskip 0.1cm

{\bf Theorem.} (Solovay)  In primitive recursive arithmetic {\bf PRA} the following are equivalent:\footnote{It is not certain that the equiconsistency of {\bf BSK} + MP with {\bf BSK} can be proved in {\bf PRA}, but cf. \cite{jrm2019}.}
\begin{enumerate}
\item[(a)]{Con({\bf BI}).}
\item[(b)]{Con({\bf BSK} + MP).}
\item[(c)]{Con({\bf FIM} + MP).}
\end{enumerate}
A careful examination of his proof reveals the

\vskip 0.1cm

{\bf Corollary.}  In primitive recursive arithmetic {\bf PRA} the following are equivalent:
\begin{enumerate}
\item[(a)]{Con({\bf BI}).}
\item[(b)]{Con({\bf IRA} + BI$_1$ + DNS$_1$).}
\item[(c)]{Con({\bf FIM} + DNS$_1$).}
\item[(d)]{Con({\bf FIM} + MP).}
\end{enumerate}

Here DNS$_1$ is the double-negation-shift principle
\[\mathrm{\forall \alpha \neg \neg \exists x \rho(\overline{\alpha}(x)) = 0 \rightarrow \neg \neg \forall \alpha \exists x \rho(\overline{\alpha}(x)) = 0},\]
which is weaker than MP over {\bf IRA}.  DNS$_1$ is adequate for the negative interpretation of BI$_1$ over {\bf IRA} + BI$_1$ (cf. \cite{jrm2019ms}), and {\bf IRA} + BI$_1$ is adequate for Kleene's consistency proof of {\bf FIM} as Troelstra and Solovay independently observed.  Like MP, DNS$_1$ is self-realizing over {\bf IRA}.

\bibliographystyle{plain}
\bibliography{solovaynote}

\begin{thebibliography}{1}

\bibitem{kl1969}
S.~C. Kleene.
\newblock {\em Formalized recursive functionals and formalized realizability}.
\newblock Number~89 in Memoirs. Amer. Math. Soc., 1969.

\bibitem{klve1965}
S.~C. Kleene and R.~E. Vesley.
\newblock {\em The {F}oundations of {I}ntuitionistic {M}athematics,
  {E}specially in {R}elation to {R}ecursive {F}unctions}.
\newblock North Holland, 1965.

\bibitem{jrm2003}
J.~R. Moschovakis.
\newblock Classical and constructive hierarchies in extended intuitionistic
  analysis.
\newblock {\em Jour. Symb. Logic}, 68:1015--1043, 2003.

\bibitem{jrm2019ms}
J.~R. Moschovakis.
\newblock Calibrating the negative interpretation.
\newblock Extended abstract for 12th Panhellenic Logic Symposium, 2019.

\bibitem{jrm2019}
J.~R. Moschovakis.
\newblock Markov{'}s principle and subsystems of intuitionistic analysis.
\newblock {\em Jour. Symb. Logic}, 84:870--876, 2019.

\bibitem{si2009}
S.~G. Simpson.
\newblock {\em Subsystems of {S}econd {O}rder {A}rithmetic}.
\newblock Perspectives in logic. ASL, Cambridge University Press, second
  edition, 2009.

\bibitem{tr1973}
A.~S. Troelstra.
\newblock Intuitionistic formal systems.
\newblock In A.~S. Troelstra, editor, {\em Metamathematical {I}nvestigation of
  {I}ntuitionistic {A}rithmetic and {A}nalysis}, Lecture Notes in Math.
  Springer-Verlag, 1973.

\bibitem{gvfms}
G.~Vafeiadou.
\newblock A comparison of minimal systems for constructive analysis.
\newblock arXiv:1808.000383.

\end{thebibliography}

\end{document}